\documentclass[12pt]{amsart}
\usepackage{amscd, amsfonts, amssymb}
\usepackage{amscd}
\usepackage{epsfig}
\usepackage{amssymb}
\usepackage[mathscr]{eucal}
\usepackage{verbatim}
\usepackage{fullpage}
\usepackage{latexsym}
\usepackage{lscape}

\newcommand\inverse{{^{-1}}}

\newcommand\ra{\rightarrow}

\DeclareMathOperator{\Aut}{Aut}
\DeclareMathOperator{\Char}{char}

\DeclareMathOperator{\graph}{graph}

\DeclareMathOperator{\GL}{GL}

\DeclareMathOperator{\SL}{SL}

\DeclareMathOperator{\SO}{SO}
\DeclareMathOperator{\SP}{Sp}

\DeclareMathOperator{\Ker}{Ker}

\DeclareMathOperator{\rk}{rk}

\numberwithin{equation}{section}

\newtheorem{thm}[equation]{Theorem}

\newtheorem{lem}[equation]{Lemma}
\newtheorem{cor}[equation]{Corollary}
\newtheorem{prop}[equation]{Proposition}

\newtheorem{qn}[equation]{Question}   % BM

\theoremstyle{definition}
\newtheorem{defn}[equation]{Definition}
\newtheorem{exmp}[equation]{Example}
\theoremstyle{remark}
\newtheorem{rem}[equation]{Remark}
\theoremstyle{remark}


\title[Complete Reducibility and Commuting Subgroups]
{Complete Reducibility and Commuting Subgroups}
\author[M.\  Bate]{Michael Bate}
\address%[M.\  Bate]
{Christ Church College, Oxford University, Oxford, OX1 1DP, UK}
\email{bate@maths.ox.ac.uk}
\author[B.\ Martin]{Benjamin Martin}
\address%[B.\ Martin]
{Mathematics and Statistics Department,
University of Canterbury,
Private Bag 4800,
Christchurch 1,
New Zealand}
\email{B.Martin@math.canterbury.ac.nz}
\author[G. R\"ohrle]{Gerhard R\"ohrle}
\address%[G.~R\"{o}hrle]
{Fakult\"at f\"ur Mathematik,
Ruhr-Universit\"at Bochum,
  D-44780 Bochum, Germany}
\email{gerhard.roehrle@rub.de}

\begin{document}

\begin{abstract}
Let $G$ be a reductive linear algebraic group
over an algebraically closed field of characteristic $p \ge 0$.
We study J-P.\ Serre's notion of
$G$-complete reducibility for subgroups of $G$.
Specifically, for a subgroup $H$ and a normal subgroup $N$ of $H$,
we look at the relationship between $G$-complete reducibility
of $N$ and of $H$, and show that these properties are equivalent
if $H/N$ is linearly reductive, generalizing a result of Serre.
We also study the case when $H=MN$ with $M$ a $G$-completely reducible
subgroup of $G$ which normalizes $N$.
In our principal result we show that if $G$ is connected,
$N$ and $M$ are connected commuting $G$-completely reducible subgroups of $G$,
and $p$ is good for $G$,
then $H = MN$ is also $G$-completely reducible.
\end{abstract}

\maketitle

\section{Introduction}
\label{sec:intro}
Let $G$ be a reductive algebraic group defined over an
algebraically closed field of characteristic $p \ge 0$
and suppose $H$ is a subgroup of $G$.
Following Serre \cite{serre1},
we say that $H$ is \emph{$G$-completely reducible} (or simply \emph{$G$-cr})
if whenever $H$ is contained in a parabolic subgroup $P$ of $G$,
then $H$ is contained in a Levi subgroup
of $P$ (in fact, we slightly extend Serre's definition to cover
the case when $G$ is not connected,
see Section \ref{sec:prelims} for precise details).
The notion of $G$-complete reducibility was introduced by
J-P.\ Serre as a way of generalizing the notion
of complete reducibility (semisimplicity) from representation theory; indeed,
when $G = \GL(V)$, a subgroup $H$ is $G$-completely reducible if and only
if $V$ is a semisimple module for $H$ \cite{serre1}.

This paper builds on the
following result \cite[Thm.\ 3.10]{BMR},
see also \cite[Thm.\ 2]{martin2}:
\begin{thm}
\label{thm:BMR3.10}
Let $H$ be a closed subgroup of $G$ with closed normal subgroup $N$.
If $H$ is $G$-completely reducible, then so is $N$.
\end{thm}
In the case $G=\GL(V)$, this statement
reduces to Clifford's Theorem from
representation theory.
Theorem \ref{thm:BMR3.10} answers a question raised by J-P.\ Serre
\cite[p.\ 24]{serre1}, who also
provides a partial converse \cite[Property 5]{serre1} under some restrictions
on the quotient $H/N$;
other partial converses are provided by
\cite[Thm.\ 3.14, Cor.\ 3.16]{BMR}, for example.
Note the converse is not true in full generality: e.g., take
$N = \{1\}$ and $H$ a non-trivial unipotent subgroup of $G^0$
(see Remark \ref{rem:gcr=>red}).

In this paper we investigate partial converses to Theorem \ref{thm:BMR3.10}
under various restrictions.
In Section \ref{sec:clifford} we show that if the quotient group
$H/N$ is linearly reductive, then $H$ is $G$-completely reducible if and
only if $N$ is
$G$-completely reducible (see Corollary \ref{cor:linred});
when $H/N$ is finite, our result gives Serre's
converse as a special case.
We consider the following question.

\begin{qn}
\label{qn:MnormalisesN}
Let $H$, $N$ be subgroups of $G$ with $N$ normal in $H$.  Let $M$ be any
subgroup of $H$ such that $MN=H$.  Is it true that $H$ is
$G$-completely reducible if and only if $M$
and $N$ are $G$-completely reducible?
\end{qn}

If $H$ is $G$-cr and
$M$ is also normal in $H$, then
$M$ and $N$ are $G$-cr, by Theorem \ref{thm:BMR3.10}.
When $M$ is not assumed to be normal in $H$, it is easy
to construct examples where $H$ is $G$-cr but $M$ is not:
we can just take $G = H = N$ and $M$ to be a non-$G$-cr subgroup of $G$.
There are examples even when $M$ is a complement to $N$ in $H$
(Examples \ref{exmp:S3} and \ref{exmp:badconn});
the problem here is that $N$ can fail to normalize $M$.

Now consider the other implication of Question \ref{qn:MnormalisesN}.
Unfortunately, even in the best possible case when $M$ and $N$ are
connected disjoint
commuting subgroups of $G$, so that $M$ is a complement to $N$ and
$M$ is normal in $H$, the answer is no (see Example \ref{exmp:liebeck}).
However, this is a low characteristic phenomenon,
as we show in our main result:

\begin{thm}
\label{thm:commuting}
Suppose that $G$ is connected and that $p$ is good for $G$ or $p>3$.
Let $A$ and $B$ be commuting connected $G$-completely reducible subgroups of
$G$.
Then $AB$ is $G$-completely reducible.
\end{thm}

Theorem \ref{thm:commuting} was first proved in
\cite[Prop.\ 40]{mcninch1} under the assumption that
$p> 2h - 2$, where $h$ is the Coxeter number of $G$;
this bound stems from Serre's result \cite[Cor.\ 5.5]{serre2}
used in the proof of \cite[Prop.\ 40]{mcninch1}.
It follows from Theorem \ref{thm:commuting} that
\cite[Thm.\ 2]{mcninch1} in fact holds for $p$ good.

This paper is organized as follows.  In Section \ref{sec:prelims}
we recall some background material, mostly taken from \cite{BMR}
and \cite{liebeckseitz}.  In Section \ref{sec:clifford} we prove
some results on $G$-complete reducibility relevant to
Theorem \ref{thm:BMR3.10}.  Section \ref{sec:commuting}
contains the proof of Theorem \ref{thm:commuting}.
Here we rely heavily on the exhaustive work of Liebeck and Seitz
\cite{liebeckseitz}, which is based on intricate case-by-case arguments;
we blend further case-by-case arguments with the general results from
the previous sections.  In Section \ref{sec:examples} we consider some
counterexamples to the statement of Theorem \ref{thm:commuting} with
the assumption of connectedness or the hypothesis on the
characteristic removed.

\section{Notation and Preliminaries}
\label{sec:prelims}

\subsection{Basic Notation}

Throughout, we work over an algebraically closed field $k$ of characteristic
$p \geq 0$; we let $k^*$ denote the multiplicative group of $k$.
By a subgroup of an algebraic group we mean a closed subgroup and by a
homomorphism of algebraic groups we mean a group homomorphism which
is also a morphism of algebraic varieties.
Let $H$ be a linear algebraic group.
We denote by  $\overline {\langle S\rangle}$
the algebraic subgroup of $H$ generated by a subset $S$ of $H$.
We let $Z(H)$ denote the centre of $H$ and $H^0$ the connected component of
$H$ that contains $1$.
If $K$ is a subgroup of $H$, then $C_H(K)$ is the centralizer of $K$ in $H$
and $N_H(K)$ is the normalizer of $K$ in $H$.
For elements $a$ and $b$ of $H$, we denote the \emph{commutator}
$aba\inverse b\inverse$ by $[a,b]$; similarly, for subgroups
$A$ and $B$ of $H$, we denote the \emph{commutator subgroup}
$\overline{\langle [a,b] \mid a\in A, b\in B\rangle}$ by $[A,B]$.
%We let $\Lie H$ denote the Lie algebra of $H$;
% recall that $\Lie H = \Lie H^0$.

If $f : H_1 \to H_2$ is a homomorphism of algebraic groups, then we say $f$
is \emph{non-degenerate} if $(\Ker f)^0$ is a torus; in particular, an isogeny
(an epimorphism with finite kernel) is non-degenerate.
We denote the \emph{rank} of $H$ by $\rk H$.
For the set of cocharacters (one-parameter subgroups) of $H$ we write $Y(H)$;
the elements of $Y(H)$ are the homomorphisms from $k^*$ to $H$.

The \emph{unipotent radical} of $H$ is denoted $R_u(H)$; it is the maximal
connected normal unipotent subgroup of $H$.
The algebraic group $H$ is called \emph{reductive} if $R_u(H) = \{1\}$;
note that we do not insist that a reductive group is connected.
In particular, $H$ is reductive if it is simple as an algebraic group
($H$ is said to be \emph{simple} if $H$ is connected and the only  
proper normal subgroups of $H$ are finite).
For a connected reductive group $H$, we can write $[H,H] = H_1\cdots  
H_t$, where
the $H_i$ are commuting simple subgroups of $H$; we call these subgroups
the \emph{simple factors} of $H$.
The algebraic group $H$ is called \emph{linearly reductive} if all rational
representations of $H$ are semisimple.
It is known that if $p = 0$, then $H$ is linearly reductive
if and only if $H$ is reductive;
in contrast, M.\ Nagata has shown (\cite[\S4, Thm.\ 2]{nagata})
that if $p > 0$,
then $H$ is linearly reductive if and only if $H^0$ is a torus and $H/H^0$ has
order coprime to $p$.
%Equivalently, this holds if and only if $H$ consists entirely of semisimple
%elements.

Throughout the paper $G$ denotes a reductive algebraic group, possibly
non-connected.
%We sometimes use the notation $\gg$ for $\Lie G$.
A subgroup of $G$ normalized by some maximal torus $T$
of $G$ is called a \emph{regular} subgroup of $G$ (connected semisimple
regular subgroups of connected reductive groups are often also referred to as
\emph{subsystem subgroups}, e.g.,\ see \cite{liebeckseitz}).

Let $\Psi = \Psi(G,T)$ denote the set of roots of $G$
with respect to a maximal torus $T$.
Fix a Borel subgroup $B$ of $G$ containing $T$ and let
$\Sigma = \Sigma(G, T)$
be the set of simple roots of $\Psi$ defined by $B$. Then
$\Psi^+ = \Psi(B, T)$ is the set of positive roots of $G$.
For $\beta \in \Psi^+$ write
$\beta = \sum_{\alpha \in \Sigma} c_{\alpha\beta} \alpha$
with $c_{\alpha\beta} \in \mathbb N_0$.
A prime $p$ is said to be \emph{good} for $G$
if it does not divide any non-zero $c_{\alpha\beta}$, and {\em bad} otherwise.
A prime $p$ is good for $G$ if and only if it is good for every simple factor
of $G$ \cite{SS}; the bad primes for the simple groups are $2$ for all groups
except type $A_n$, $3$ for the exceptional groups and $5$ for type $E_8$.

Suppose $G$ acts on a variety $V$ and let $v \in V$.
Then for each cocharacter $\lambda \in Y(G)$, we can
define a morphism of varieties
$\phi_{v,\lambda}:k^* \to V$ via the formula
$\phi_{v,\lambda}(x) = \lambda(x)\cdot v$.
If this morphism extends to a morphism
$\overline\phi_{v,\lambda}:k \to V$, then
we say that $\underset{x\to 0}{\lim}\, \lambda(x) \cdot v$ exists,
and set this limit
equal to $\overline\phi_{v,\lambda}(0)$; note that such an extension,
if it exists, is necessarily unique.
This procedure is especially important for us when $V=G$
and $G$ acts by conjugation.
 
\subsection{$G$-Complete Reducibility}
\label{subsec:noncon}
In \cite[Sec.\ 6]{BMR}, Serre's original
notion of $G$-complete reducibility is extended
to include the case where $G$ is reductive but not necessarily connected
(so that $G^0$ is a connected reductive group).
The crucial ingredient of this extension is the introduction of so-called
\emph{Richardson parabolic subgroups} (\emph{R-parabolic subgroups}) of a
%non-connected
reductive group $G$.
We briefly recall the main definitions and results;
for more details and further results,
the reader is referred to \cite[Sec.\ 6]{BMR}.

\begin{defn}
\label{defn:rpars}
For each cocharacter $\lambda \in Y(G)$, let
$P_\lambda = \{ g\in G \mid \underset{x \to 0}{\lim}\,
\lambda(x) g \lambda(x)\inverse \textrm{ exists} \}$.
Recall that a subgroup $P$ of $G$ is \emph{parabolic}
if $G/P$ is a complete variety.
The subgroup $P_\lambda$ is parabolic in this sense, but the converse
is not true: e.g.,\ if $G$ is finite,
then every subgroup is parabolic, but the only
subgroup of $G$ of the form $P_\lambda$ is $G$ itself.
If we define
$L_\lambda = \{g \in G \mid \underset{x \to 0}{\lim}\,
\lambda(x) g \lambda(x)\inverse = g\}$,
then $P_\lambda = L_\lambda \ltimes R_u(P_\lambda)$,
and we also have
$R_u(P_\lambda) = \{g \in G \mid \underset{x \to 0}{\lim}\,
\lambda(x) g \lambda(x)\inverse = 1\}$.
The subgroups $P_\lambda$ for $\lambda \in Y(G)$
are called the \emph{R-parabolic subgroups}
of $G$.
Given an R-parabolic subgroup $P$,
an \emph{R-Levi subgroup} of $P$ is any subgroup $L_\lambda$
such that $\lambda \in Y(G)$ and $P=P_\lambda$.
If $G$ is connected, then the R-parabolic subgroups
(resp.\ R-Levi subgroups of R-parabolic subgroups)
of $G$ are exactly the parabolic subgroups
(resp.\ Levi subgroups of parabolic subgroups) of $G$;
indeed, most of the theory of parabolic subgroups and
Levi subgroups of connected reductive groups
carries over to R-parabolic and R-Levi subgroups of
arbitrary reductive groups.
In particular, all R-Levi subgroups of an R-parabolic
subgroup $P$ are conjugate under the action
of $R_u(P)$.
\end{defn}

We are often interested in reductive subgroups of reductive groups.
If $H$ is a subgroup of $G$, then there is an obvious
inclusion $Y(H) \subseteq Y(G)$ of the sets of cocharacters.
When $H$ is reductive and $\lambda \in Y(H)$, there is then an R-parabolic
subgroup of $H$ associated to $\lambda$, as well as an R-parabolic
subgroup of $G$.
In order to distinguish between R-parabolic subgroups associated to  
different subgroups of $G$,
we use the notation
$P_\lambda(H)$, $L_\lambda(H)$, etc.\ where necessary, but we write
$P_\lambda$ for $P_\lambda(G)$ and $L_\lambda$ for  $L_\lambda(G)$.
Note that $P_\lambda(H) = P_\lambda \cap H$,
$L_\lambda(H) = L_\lambda \cap H$
and $R_u(P_\lambda(H)) = R_u(P_\lambda) \cap H$.

We now have the machinery in place to define exactly what we mean by a
$G$-completely reducible subgroup  in this more general setting;
note that, by the remarks above,
the following definition coincides with Serre's notion for connected $G$.

\begin{defn}
\label{defn:gcr}
  Suppose $H$ is a subgroup of $G$.
We say $H$ is \emph{$G$-completely reducible} ($G$-cr for short)
if whenever $H$ is contained
in an R-parabolic subgroup of $G$, then there exists an R-Levi subgroup
$L$ of $P$ with $H \subseteq L$.
Equivalently, $H$ is $G$-completely reducible if
whenever $H \subseteq P_\lambda$ for
some $\lambda \in Y(G)$, there exists $\mu \in Y(G)$ such that
$P_\lambda = P_\mu$ and
$H \subseteq L_\mu$.  Since all R-Levi subgroups of an R-parabolic  
subgroup $P$ are $R_u(P)$-conjugate, we have another formulation: $H$  
is $G$-completely reducible if
whenever $H \subseteq P_\lambda$ for
some $\lambda \in Y(G)$, there exists $u\in R_u(P_\lambda)$ such that
$H \subseteq uL_\lambda u^{-1}$.

Note that, for $G=\GL(V)$, a subgroup $H$ is $G$-completely reducible
if and only if $V$ is a semisimple $H$-module.
\end{defn}

\begin{rem}
\label{rem:gcr=>red}
If $G$ is connected, then any $G$-completely reducible
subgroup $H$ is reductive: in fact $H$ cannot normalize any non-trivial
unipotent subgroup of $G^0$.
This was shown by Serre in \cite[\S 4.1]{serre2}.
Serre's argument gives the same result for non-connected $G$,
cf.\ \cite[Prop.\ 5.4(b)]{martin1}.
\end{rem}

Sometimes we come across subgroups of $G$
which are not contained in any R-parabolic subgroup of $G$;
these subgroups are trivially $G$-completely reducible.
Following Serre again, we call these subgroups
\emph{$G$-irreducible} ($G$-ir).
Note that, for $G=\GL(V)$, a subgroup $H$ is $G$-irreducible if and  
only if $V$
is an irreducible $H$-module.

\smallskip

We recall some results on $G$-complete reducibility (mainly from \cite{BMR})
which we require in the sequel.
The reader should note that many of
the results in \cite{BMR} were proved first for connected
groups and then generalized to the non-connected case;
we give separate references where appropriate from \cite{BMR}
for the connected and non-connected versions of each result (see the  
paragraph immediately preceding \cite[Sec.\ 6.2]{BMR}).

Our first result (see \cite[Lem.\ 11.24]{Jantzen} for a proof),
together with Remark \ref{rem:gcr=>red},
shows that in characteristic zero,
a subgroup of $G$ is $G$-completely reducible if and only if it is
reductive; thus our results only have independent interest in
positive characteristic. This equivalence
is not true in positive characteristic.
For we can take a reductive but not linearly reductive group $H$ and map it
into some $\GL(V)$ in a non-semisimple way; the image of $H$
is then reductive but not $\GL(V)$-cr.

\begin{lem}
\label{lem:linred=>gcr}
If $H$ is a linearly reductive
subgroup of $G$, then $H$ is $G$-completely reducible.
\end{lem}

Regular subgroups play an important r\^ole in many of our proofs;
the following result \cite[Prop.\ 3.20]{BMR} shows that these  
subgroups are $G$-completely reducible.

\begin{lem}
\label{lem:reg=>gcr}
  If $H$ is a regular reductive subgroup of $G$, then $H$ is  
$G$-completely reducible.
\end{lem}

We need a useful preliminary result.

\begin{lem}
\label{lem:almostsplit}
  Let $N$ be a normal subgroup of $G$.  Then there exists a subgroup
$M$ of $G$ such that $MN=G$, $M\cap N$ is a finite normal subgroup
of $M$, $M^0\cap N^0$ is central in both $M^0$ and $N^0$, and $M^0$
commutes with $N$.
\end{lem}

\begin{proof}
  The existence of $M$ and all of its properties except the last
follow from \cite[Lem.\ 6.14]{BMR} and its proof.  To show
that $M^0$ commutes with $N$, we observe that for any $x\in N$, the
connected set $\{[x,y] \mid y \in M^0\}$ is contained in the finite
set $M\cap N$ and hence must be trivial.
\end{proof}

The next result (see \cite[Lem.\ 2.12, Sec.\ 6.2]{BMR}) shows that the  
concept of complete reducibility behaves well with respect to certain  
homomorphisms.
Observe that Lemma \ref{lem:bmr2.12}(ii)(b) applies in particular to  
isogenies;
we use this special case frequently in the sequel.

\begin{lem}
\label{lem:bmr2.12}
Let $G_1$ and $G_2$ be reductive groups.
\begin{itemize}
\item[(i)] Let $H$ be a  subgroup of $G_1 \times G_2$.
Let $\pi_i:G_1 \times G_2 \to G_i$ be the canonical projection
for $i=1,2$.
Then $H$ is $(G_1\times G_2)$-completely reducible if and only if  
$\pi_i(H)$ is $G_i$-completely reducible
for $i=1,2$.
\item[(ii)] Let $f:G_1 \to G_2$ be an epimorphism. Let $H_1$ (resp.\ $H_2$)
be a subgroup of $G_1$ (resp.\ $G_2$).
\begin{itemize}
\item[(a)] If $H_1$ is $G_1$-completely reducible,
then $f(H_1)$ is $G_2$-completely reducible.
\item[(b)] If $f$ is non-degenerate, then $H_1$ is
$G_1$-completely reducible if and only if
$f(H_1)$ is $G_2$-completely reducible, and
$H_2$ is $G_2$-completely reducible if and only if
$f\inverse(H_2)$ is $G_1$-completely reducible.
\end{itemize}
\end{itemize}
\end{lem}

In particular, let $G=G_1\times G_2$ and consider the normal
subgroup $G_1$ of $G$.  If $H$ is a subgroup of $G_1$, then $H$
is $G$-completely reducible if and only if $H$ is $G_1$-completely reducible.
We extend this result to cover the case of a normal subgroup which is not
necessarily a direct factor of $G$.

\begin{prop}
\label{prop:10}
Let $K\subseteq N\subseteq G$ be subgroups of $G$ with $N$ normal in $G$.
Then $K$ is $G$-completely reducible if and only if $K$ is $N$-completely
reducible.
\end{prop}

\begin{proof}
Let $M$ be as in Lemma \ref{lem:almostsplit}.
Suppose $K$ is $G$-cr.  Let $\lambda\in Y(N)$ such that
$K\subseteq P_\lambda(N)$.  As $K$ is $G$-cr, there exists
$u\in R_u(P_\lambda)$ such that $K\subseteq uL_\lambda u^{-1}$.
Now $R_u(P_\lambda) \subseteq G^0$, so we can write $u = u_1u_2$ with
$u_1\in N^0$ and $u_2\in M^0$.  Since $\lambda$ centralizes
$M^0$, we see that
$u_1':= \underset{x\ra 0}{\lim}\, \lambda(x)u_1\lambda(x)^{-1}$ must
exist and equal $u_2^{-1}$.  Then $u_2\in M^0\cap N^0$, which consists of
semisimple elements, so $u_2 = 1 = u_1'$ and $u_1\in R_u(P_\lambda(N))$.
Thus $K\subseteq u_1L_\lambda(N)u_1^{-1}$, an R-Levi subgroup of
$P_\lambda(N)$.  This shows that $K$ is $N$-cr.

Conversely, suppose $K$ is $N$-cr.  Let $\lambda\in Y(G)$
such that $K\subseteq P_\lambda$.
Since $M^0$ commutes with $N$, we can write $\lambda= \sigma+\tau$,
where $\sigma\in Y(M)$ and $\tau\in Y(N)$.
Since $\sigma$ centralizes $N$, $\sigma$ centralizes $K$, so
$K\subseteq P_\tau\cap N= P_\tau(N)$.  As $K$ is $N$-cr,
there exists $u\in R_u(P_\tau(N)) = R_u(P_\tau)\cap N = R_u(P_\lambda)  
\cap N$ such that
$uKu^{-1}\subseteq L_\tau(N)= L_\lambda\cap N$.  Thus $K$ is $G$-cr,
as required.
\end{proof}

The following three results allow us to relate $G$-complete reducibility
and $H$-complete reducibility for subgroups of some important  
subgroups $H$ of $G$.
The first result (\cite[Cor. 3.21]{BMR})
makes sense because, if $S$ is a linearly reductive group acting on $G$
by automorphisms, then $C_G(S)$ is reductive \cite[Prop.\ 10.1.5]{rich}.

\begin{prop}
\label{prop:linred}
Let $S$ be a linearly reductive group acting on $G$ by automorphisms  
and let $H= C_G(S)^0$.
Suppose $K$ is a  subgroup of $H$.
Then $K$ is $H$-completely reducible
if and only if $K$ is $G$-completely reducible.
\end{prop}

This is of particular use to us when $G$ is simple
and $S$
is the finite group generated by a graph automorphism of $G$.
The next result is a corollary of Proposition \ref{prop:linred}
when $G$ is connected; the point is that any R-Levi subgroup of
an R-parabolic subgroup of $G$ is the centralizer of a torus.  The extension
to non-connected groups is not hard (see \cite[Cor. 3.22, Sec.\ 6.3]{BMR}).

\begin{cor}
\label{cor:levi}
Let $L$ be an R-Levi subgroup
of some R-parabolic subgroup of $G$.
Then a subgroup $K$ of $L$ is $G$-completely reducible if and only if  
it is $L$-completely reducible.
\end{cor}

The next result is \cite[Thm.\ 3.26]{BMR} when $G$ is connected. The  
proof given there still applies for non-connected $G$ by Proposition  
\ref{prop:finiteindex} below, since $H$ is assumed to be connected.

\begin{prop}
\label{prop:regular}
Suppose that $p$ is good for $G^0$.
Let $H$ be a regular connected reductive subgroup of $G$ and suppose
$K$ is a subgroup of $H$.
Then $K$ is $H$-completely reducible if and only if
$K$ is $G$-completely reducible.
\end{prop}

The following result is a generalization of \cite[Lem.\ 6.12(i)]{BMR}.
Note that a finite-index subgroup of a reductive group is reductive.

\begin{prop}
\label{prop:finiteindex}
Let $K\subseteq H$ be subgroups of $G$,
with $H$ of finite index in $G$.
Then $K$ is $H$-completely reducible if and only if
$K$ is $G$-completely reducible.
\end{prop}

\begin{proof}
  Suppose $K$ is $G$-cr.  Let $P$ be an R-parabolic
subgroup of $H$ with $K\subseteq P$.  We can write
$P=P_\lambda(H)$ for some $\lambda\in Y(H)$.  Then
$K\subseteq P_\lambda$.  Since $K$ is $G$-cr, there
exists $u\in R_u(P_\lambda)$ such that $u K u^{-1}\subseteq L_\lambda$.
We have $R_u(P_\lambda(H)) = H\cap R_u(P_\lambda) = H^0\cap  
R_u(P_\lambda) = G^0\cap R_u(P_\lambda) = R_u(P_\lambda)$, so
$u$ belongs to $R_u(P_\lambda(H))$ and
$u K u^{-1}\subseteq H\cap L_\lambda = L_\lambda(H)$.  Thus $K$ is $H$-cr.
The argument in the other direction is similar.
\end{proof}

We obtain a corollary which strengthens Proposition \ref{prop:10}.

\begin{cor}
\label{cor:finiteconjs}
  Let $N$ be a subgroup of $G$ such that the number of $G$-conjugates
of $N$ is finite, and let $K$ be a subgroup of $N$.
Then $K$ is $G$-completely reducible if and only if
$K$ is $N$-completely reducible.
\end{cor}

\begin{proof}
The subgroup $G_1:=N_G(N)$ has finite index in $G$.
Then $G_1$ is reductive, so its normal subgroup
$N$ is also reductive.  The result now follows from Propositions
\ref{prop:10} and \ref{prop:finiteindex}.
\end{proof}

\begin{rem}
  Let $K$, $H$ be subgroups of $G$ with $H$ of finite index in $G$.
It need not be true that if $K$ is $G$-completely reducible
then $K\cap H$ is $H$-completely reducible.  For example, let
$p=2$, let $\phi$ be an irreducible embedding of the symmetric
group $S_3$ in $\SL_2$, let $G = S_3\times \SL_2$ and let
$K=\{g\phi(g) \mid g\in S_3\}$.  Let $C_2$ be
a cyclic subgroup of $S_3$ of order 2 and let
$H=C_2\times \SL_2\subseteq G$.  It is easily
checked that $K$ is $G$-ir but $K\cap H$ is not $H$-cr.
\end{rem}

Some of the proofs in Section \ref{sec:commuting}
rely on the monograph \cite{liebeckseitz} of M.\ Liebeck and G.\ Seitz.
We recall some notation and results from this paper.
For the rest of this section, assume $G$ is simple (hence connected) and
of exceptional type.  For convenience, we take $p$ to be $\infty$  
rather than zero if $k$ has characteristic zero.
For certain simple subgroups $X$
of $G$, Liebeck and Seitz provide a positive integer
$N(X,G)$ given in the table at the top of page 2 of
\cite{liebeckseitz}. We reproduce
this table for ease of reference (Table \ref{t:nxg}).

\begin{table}[h]
\begin{tabular}{r|rrrrr}
         & $G=E_8$   & $E_7$     & $E_6$     & $F_4$     & $G_2$\\ \hline
$X= A_1$& $7$       & $7$       & $5$       & $3$       & $3$  \\
$A_2$   & $5$       & $5$       & $3$       & $3$       &      \\
$B_2$   & $5$       & $3$       & $3$       & $2$       & \\
$G_2$   & $7$       & $7$       & $3$       & $2$       & \\
$A_3$   & $2$       & $2$       & $2$       &           & \\
$B_3$   & $2$       & $2$       & $2$       & $2$       & \\
$C_3$   & $3$       & $2$       & $2$       & $2$       & \\
$B_4, C_4, D_4$ & $2$ & $2$     & $2$       &           &
\end{tabular}
\bigskip
\caption{The bounds $N(X,G)$ from \cite{liebeckseitz}.}
\label{t:nxg}
\end{table}
For example, if $X$ has type $A_3$ and $G$ has type $E_6$, then $N(X,G) = 2$.
If the pair $(X,G)$ is not in the table, then set $N(X,G) = 1$.
More generally, if $X$ is a connected reductive subgroup of
$G$, and $X_1, \ldots, X_t$ are the simple factors of $X$,
then define
$$
N(X,G) = \max(N(X_i,G) \mid 1 \leq i \leq t),
$$
where we take $N(X,G) = 1$ if $X$ is a torus.
With this definition, we can restate \cite[Thm.\ 1, Thm.\  
3.8]{liebeckseitz} in
our language.

\begin{thm}
\label{thm:ls1}
Let $X$ be a  connected reductive subgroup of $G$ and assume that
$p > N(X,G)$. Then $X$ is $G$-completely reducible.
\end{thm}

\begin{thm}
\label{thm:ls3.8}
Let $C$ be a simple subgroup of $G$ such that $C$ is
of classical type.
Suppose that $X$ is a  connected reductive subgroup of $C$
and that $p>N(X,G)$.
When $(X,p) = (B_l,2)$ or $(C_l,2)$, assume that $C \neq B_r$ or $C_r$.
Then $X$ is $C$-completely reducible.

\end{thm}

Note that whether a given $X$ is $G$-completely reducible depends only on the
Dynkin type of  $G$ and not on its isogeny class (Lemma  
\ref{lem:bmr2.12}(ii)(b)),
so there is no harm in labelling $G$ only by its Dynkin type.

\begin{rem}
\label{rem:useoftables}
Tables 8.1--8.5 of \cite{liebeckseitz}
give the connected simple subgroups $X$ of an exceptional
group $G$, the connected centralizers $C_G(X)^0$, and the minimal
connected semisimple regular subgroups of $G$ that contain $X$,
under the hypothesis that $p>N(X,G)$.  Note that these tables
give $X$ up to $\Aut(G)$-conjugacy.  We give an example
to point out one consequence of this.  Suppose $p>2$.  A group
of type $D_6$ contains a subgroup of type $B_5$ as the connected
centralizer of a graph automorphism of order $2$.
A group of type $E_7$ contains a group of type $D_6$ as a Levi subgroup.
 From Table 8.2 of \cite{liebeckseitz},
we see that there is exactly one copy up to $\Aut(G)$-conjugacy
of a group $X$ of type $B_5$ in a group $G$ of type $E_7$,
so it must be the subgroup we
already know: that is, it must be the connected centralizer of a
graph automorphism of order $2$ of the $D_6$-subgroup of $E_7$.
We use arguments of this kind repeatedly in the proof of
Theorem \ref{thm:commuting}.
\end{rem}

\section{Converses to Theorem \ref{thm:BMR3.10}}
\label{sec:clifford}

In this section we give several results providing conditions
under which the converse of Theorem \ref{thm:BMR3.10} holds, culminating
in Corollary \ref{cor:linred}, which generalizes a result of Serre.
We also prove some technical results which prepare the ground for
the proof of our main result in Section \ref{sec:commuting}.
We begin with a refinement of \cite[Prop.~3.19]{BMR}.

\begin{thm}
\label{thm:HKM}
Suppose $K\subseteq H \subseteq M$ are subgroups of $G$ such that
$M$ is reductive.
\begin{itemize}
  \item[(a)] Suppose $M$ contains a maximal torus of $C_G(K)$.
  \begin{itemize}
   \item[(i)] If $H$ is $G$-completely reducible,
then $M$ is $G$-completely reducible.
   \item[(ii)] If $K$ is $G$-completely reducible,
then $M$ is $G$-completely reducible.
  \end{itemize}
  \item[(b)] Suppose $C_G(K)^0 \subseteq M$.
  \begin{itemize}
   \item[(i)] If $H$ is $G$-completely reducible,
then $H$ is $M$-completely reducible.
   \item[(ii)] If $K$ is $G$-completely reducible
and $H$ is $M$-completely reducible, then $H$ is $G$-completely reducible.
  \end{itemize}
\end{itemize}
\end{thm}

\begin{proof}
(a). Note that part (ii) follows from part (i),
just taking $H$ to be $K$.  To prove part (i),
suppose $P$ is an R-parabolic subgroup of $G$ containing $M$.
Then $H$ lies in an R-Levi subgroup $L$ of $P$, since $H$ is
$G$-cr, so $K$ lies in $L$.  Let $S$ be a maximal torus
of $C_G(K)$ such that $S\subseteq M$.  Then $S$ is a maximal torus of  
$C_P(K)$, so, after conjugating
$L$ by some element of $C_P(K)$, we can assume that $Z(L)^0\subseteq  
S$.  Pick $\lambda\in Y(G)$ such that $P=P_\lambda$ and $L=L_\lambda$.  
  Then $\lambda\in Y(S)\subseteq Y(M)$.  We have $P_\lambda(M)=  
P_\lambda\cap M= M$, so $\lambda\in Y(Z(M))$, by \cite[Lem.\  
2.4]{BMR}.  Thus $M\subseteq L_\lambda=L$.  This shows that $M$ is  
$G$-cr.

(b)(i). We have $C_G(H)^0\subseteq C_G(K)^0\subseteq M$ by hypothesis.
Now \cite[Prop.~3.19]{BMR} implies that if $H$ is $G$-cr,
then $H$ is $M$-cr.

(ii). Suppose $P$ is an R-parabolic subgroup of $G$ containing $H$;
then $K\subseteq P$ also.
We are given that $K$ is $G$-cr, so there exists an R-Levi subgroup $L$ of $P$
with $K \subseteq L$.
Let $\lambda \in Y(G)$ be such that $P = P_\lambda$ and $L = L_\lambda$.
Then, since $K \subseteq L = C_G(\lambda(k^*))$, we have
$\lambda \in Y(C_G(K)^0) \subseteq Y(M)$.
Thus $P_\lambda(M)$ is an R-parabolic subgroup of $M$ containing $H$.
Since $H$ is $M$-cr, there exists
$u \in R_u(P_\lambda(M))=R_u(P_\lambda)\cap M$ such that
$H\subseteq u L_\lambda(M) u^{-1}$.  Hence $H\subseteq uL_\lambda u^{-1}$,
an R-Levi subgroup of $P$.  This shows that $H$ is $G$-cr.
\end{proof}

\begin{rem}
  (a) Note that in part (a) of Theorem \ref{thm:HKM} we only require $M$
to contain a maximal torus of $C_G(K)$,
rather than all of $C_G(K)^0$.
Part (b) fails under this weaker hypothesis, however.  For example,
we can take $K$ to be $\{1\}$ and $H$ to be a subgroup of $M$, where
$M$ contains a maximal torus of $G$.  There exist examples in which
$H$ is $M$-cr but not $G$-cr, and others in which $H$ is $G$-cr but  
not $M$-cr.
For the former, see \cite[Ex.\ 3.45]{BMR}: we take $H$ to be
$\SP_m$ embedded diagonally in $M:=\SP_m\times \SP_m$ inside
$G:=\SP_{2m}$, where $m \ge 4$ is even.
For the latter, we can take $H$ and $M$ to be
certain subgroups of a simple group of type $G_2$ in characteristic $2$;
see \cite[Prop.\ 7.17]{BMRT}.

  (b) Putting $M=H$ in Theorem \ref{thm:HKM}, we obtain the following  
result: if $K$ is a $G$-cr subgroup of $G$, then any reductive  
subgroup $H$ of $G$ containing $KC_G(K)^0$ is also $G$-cr (note that  
$H$ is $H$-cr).  This is a strengthening of \cite[Thm.\ 3.14]{BMR}.
\end{rem}

\begin{cor}
\label{cor:gcr<=>ncr}
Suppose $N \subseteq H$ are subgroups of $G$,
$N$ is normal in $H$ and $N$ is $G$-completely reducible.
Then $H$ is $G$-completely reducible if and only if $H$ is
$N_G(N)$-completely reducible.
\end{cor}

\begin{proof}
This follows from Theorem \ref{thm:HKM}(b), setting $K = N$ and $M = N_G(N)$.
Note that since $N$ is $G$-cr, $N_G(N)$ is %non-connected
reductive, by the non-connected version of \cite[Prop.\ 3.12]{BMR}.
\end{proof}

%Recall that we call a homomorphism of algebraic groups
%\emph{non-degenerate} if its connected kernel is a torus.

\begin{thm}
\label{thm:gmodhcr}
Suppose $N \subseteq H$ are subgroups of $G$ and $N$ is normal in $G$.
Then $H$ is $G$-completely reducible if and only if
$H/N$ is $G/N$-completely reducible.
\end{thm}

\begin{proof}
Let $\pi: G \to G/N$ denote the canonical homomorphism.
We first perform a series of reductions, using results from \cite{BMR}.
Let $M$ be as in Lemma \ref{lem:almostsplit}.
Since $N$ is normal in $G$, $Z(N^0)^0$ is normal in $G$;
the canonical map $G \to G/Z(N^0)^0$ is non-degenerate, so by Lemma  
\ref{lem:bmr2.12}(ii)(b)
we may assume that $Z(N^0)^0 = \{1\}$.
The map $M \ltimes N \to G$ induced by multiplication
is also non-degenerate,
so we may assume that $G = M \ltimes N$ and, therefore,
that the map $\pi: G \to G/N$ is the projection onto the first factor.

Now suppose $\lambda \in Y(M)$.
Then, since $N$ and $M^0$ commute,
$\lambda(k^*)$ centralizes $N$.
By \cite[Lem.\ 6.15(i),(ii)]{BMR}, we have
\begin{align}
\pi(P_\lambda) = P_\lambda(M), \qquad \pi\inverse(P_\lambda(M)) =
P_\lambda,
\label{eq:piP}\\
\pi(L_\lambda) = L_\lambda(M), \qquad \pi\inverse(L_\lambda(M)) =
L_\lambda.
\label{eq:piL}
\end{align}
If $P = P_\mu$ is an R-parabolic subgroup of $G$ containing $H$,
then we can write $\mu = \lambda + \nu$, where $\lambda \in Y(M)$
and $\nu \in Y(N)$.
Since $P$ contains $H$, it contains $N$; but
$N$ is normal in $G$, and hence $N$ is $G$-cr by
Theorem \ref{thm:BMR3.10}, which
means that $N$ is contained in some R-Levi subgroup of $P$.
Moreover, normality of $N$ in $P$ implies $N$ is contained in \emph{every}
R-Levi subgroup of $P$; in particular, $N \subseteq L_\mu$, so
$\mu(k^*)$ centralizes $N$.
We have already noted that $\lambda(k^*)$ centralizes
$N$ for any $\lambda \in Y(M)$.
Thus $\nu(k^*)$ centralizes $N$, and, since $Z(N^0)^0 = \{1\}$,
$\nu$ is trivial.
We can finally conclude that $\mu \in Y(M)$.

It is now easy to see from Equations
\eqref{eq:piP} and \eqref{eq:piL} that
$H$ is contained in an R-parabolic (resp.\ R-Levi) subgroup of $G$ if  
and only if
$H/N$ is contained in the corresponding R-parabolic (resp.\ R-Levi) subgroup
of $G/N$. Thus $H$ is $G$-cr if and only if $H/N$ is $G/N$-cr, as required.
\end{proof}

We can now prove the main result of this section.

\begin{cor}
\label{cor:linred}
Suppose that $N \subseteq H$ are subgroups of $G$ with $N$ normal in  
$H$. %Then
\begin{itemize}
\item[(i)] If $N$ is $G$-completely reducible,
then $H$ is $G$-completely reducible
if and only if $H/N$ is $N_G(N)/N$-completely reducible.
\item[(ii)] If $H/N$ is linearly reductive,
then $H$ is $G$-completely reducible if and only if
$N$ is $G$-completely reducible.
\end{itemize}
\end{cor}

\begin{proof}
(i).
By Corollary \ref{cor:gcr<=>ncr}, $H$ is $G$-cr
if and only if $H$ is $N_G(N)$-cr.
By Theorem \ref{thm:gmodhcr} applied to the inclusions $N \subseteq H  
\subseteq N_G(N)$,
it follows that $H$ is $N_G(N)$-cr if and only if
$H/N$ is $N_G(N)/N$-cr.
This proves the result.

(ii).
Suppose $H/N$ is linearly reductive.
If $N$ is $G$-cr, then $N_G(N)$ is a reductive group and
so is the quotient $N_G(N)/N$.
Thus  $H/N$ is a linearly reductive subgroup of $N_G(N)/N$.
By Lemma \ref{lem:linred=>gcr}, $H/N$ is $N_G(N)/N$-cr,
and hence $H$ is $G$-cr, by part (i).
On the other hand, if $H$ is $G$-cr, then $N$ is $G$-cr, by Theorem
\ref{thm:BMR3.10}.
\end{proof}

\begin{rem}
Note that Serre's result \cite[Property 5]{serre1} is a special case of
the reverse implication of
Corollary \ref{cor:linred}(ii); the finite linearly reductive groups
are exactly those whose orders are coprime to $p = \Char k$.
Corollary \ref{cor:linred}(ii) answers a question posed by
B.\ K\"ulshammer.
\end{rem}

Corollary \ref{cor:linred}(ii) provides a useful criterion to
ensure that $G$-complete reducibility of $H$ and of $N$ are equivalent.
However, there are many examples where $H$ and $N$ are $G$-completely
reducible, but $H/N$ is not linearly reductive.  The problem in  
general is that $G$-complete reducibility of
$N$ and $H$ depends
not only on how $N$ sits inside $H$, but also on how
$H$ sits inside $G$.
Therefore, to make more progress, one has to impose further conditions on $H$.
In the next section we consider the case when $H = MN$,
where $M$ commutes with $N$ and $M$ is also $G$-completely reducible.
We now give some results applicable to this special case, the first of which
is used many times in the proof of Theorem \ref{thm:commuting}.

\begin{prop}
\label{prop:AB<=>CGA}
Suppose
$N$ is a $G$-completely reducible subgroup of $G$ and $M$ is a  
subgroup of $C_G(N)$. Then
$M$ is $C_G(N)$-completely reducible if and only if $MN$ is  
$G$-completely reducible.
\end{prop}

\begin{proof}
By Corollaries \ref{cor:gcr<=>ncr} and \ref{cor:linred}(i),
$MN$ is $G$-cr if and only if $MN$ is $N_G(N)$-cr if and only if  
$MN/N$ is $N_G(N)/N$-cr.
Now $C_G(N)N/N$ is normal in $N_G(N)/N$,
so $MN/N$ is $N_G(N)/N$-cr if and only if $MN/N$ is $C_G(N)N/N$-cr  
(Proposition \ref{prop:10}).
Let $f\colon C_G(N)\ra C_G(N)N/N$ be the inclusion of
$C_G(N)$ in $C_G(N)N$ followed by the canonical projection from  
$C_G(N)N$ to $C_G(N)N/N$.
The connected kernel of $f$ is $(N\cap C_G(N))^0= Z(N)^0$,
which is a torus since $N$ is reductive, and $f(M)$ is $MN/N$.
It follows from Proposition \ref{lem:bmr2.12}(ii) that
$M$ is $C_G(N)$-cr if and only if $MN/N$ is $C_G(N)N/N$-cr, which  
proves the result.
\end{proof}

\begin{cor}
\label{cor:CGAcr=>Gcr}
Suppose $N$ is $G$-completely reducible and $M$ is a subgroup of $C_G(N)$.
If $M$ is $C_G(N)$-completely reducible, then $M$ is $G$-completely reducible.
\end{cor}

\begin{proof}
Under these hypotheses, $MN$ is $G$-cr, by Proposition \ref{prop:AB<=>CGA}.
But $M$ is a normal subgroup of $MN$, so $M$ is $G$-cr,
by Theorem \ref{thm:BMR3.10}.
\end{proof}

In Section \ref{sec:examples} we will show that the converse
of Corollary \ref{cor:CGAcr=>Gcr} is not true in general,
although it is true in the important special case of
Theorem \ref{thm:commuting}.

We finish this section with a result
which sometimes allows us to reduce to the case of commuting subgroups.

\begin{lem}
\label{lem:linredcomm}
Let $A$, $B$ be subgroups of $G$ such that $A$ normalizes $B$,
$[A,B]$ is linearly reductive and $[A,B]$ centralizes $A$.
Then $A$ is $G$-completely reducible if and only if $C_A(B)$
is $G$-completely reducible, and $AB$ is $G$-completely reducible
if and only if $C_A(B)B$ is $G$-completely reducible.
\end{lem}

\begin{proof}
Note that $C_A(B)$, $C_A(B)B$ are normal subgroups of $A$,
$AB$ respectively, so by Corollary \ref{cor:linred}(ii),
it is enough to show that $A/C_A(B)$ and $AB/C_A(B)B$ are linearly reductive.
By a standard DCC argument, we can choose
$b_1,\ldots,b_r\in B$ for some $r\in {\mathbb N}$ such that
$C_A(B)= C_A(\langle b_1,\ldots,b_r\rangle)$.
Define $\phi\colon A\ra [A,B]^r$ by
\[
\phi(a)= ([a\inverse,b_1],\ldots,[a\inverse, b_r]).
\]
Since $A$ commutes with $[A,B]$, we have
\[
[a\inverse,b_i][a'\inverse,b_i]
= a'\inverse[a\inverse,b_i]b_ia'b_i\inverse
= (aa')\inverse b_i (aa') b_i\inverse
= [(aa')\inverse, b_i],
\]
for each $i$ and each $a,a'\in A$.
It follows that $\phi$ is a homomorphism from $A$
to the linearly reductive group $[A,B]^r$ with kernel $C_A(B)$,
whence $A/C_A(B)$ is linearly reductive.  The obvious map from
$A/C_A(B)$ to $AB/C_A(B)B$ is surjective, so $AB/C_A(B)B$ is
also linearly reductive, as required.
\end{proof}

\section{Proof of Theorem \ref{thm:commuting}}
\label{sec:commuting}

In this section we prove our main result, Theorem \ref{thm:commuting},
via a series of reductions.
It is immediate from Lemma \ref{lem:linred=>gcr} that
Theorem \ref{thm:commuting} holds if $G$ is a torus.
The next step is to reduce to the case when $G$
is simple.  We begin with a technical definition.

\begin{defn}
\label{defn:projection}
Suppose $G$ is connected.  Write $G = G_1 \cdots G_r Z$, where
the subgroups $G_i$ are the simple factors of $G$, and $Z$ is a central
torus.
Let $\tilde{G} = G_1 \times \cdots \times G_r \times Z$,
let $\phi$ be the isogeny from $\tilde{G}$ onto $G$ induced by multiplication,
and let $\pi_i$ denote the projection of $\tilde{G}$ onto $G_i$ for each $i$.
If $X$ is a connected reductive subgroup of $G$, we call the subgroup
$\pi_i(\phi\inverse(X))$ of $G_i$ the \emph{projection of $X$ to the  
simple factor $G_i$},
and we denote it by $X_i$.
\end{defn}

\begin{lem}
\label{lem:simplefactors}
Keep the notation of Definition \ref{defn:projection}.
Then $X$ is $G$-completely reducible
if and only if $X_i$ is $G_i$-completely reducible for every $i$.
If $X$ is connected, then $X$ is $G$-completely reducible
if and only if $X_i^0$ is $G_i$-completely reducible for every $i$.
\end{lem}

\begin{proof}
By Lemma \ref{lem:bmr2.12}(ii)(b),
$X$ is $G$-cr if and only if $\phi^{-1}(X)$ is $\tilde{G}$-cr.
Let $\pi$ denote the projection of $\tilde{G}$ onto the central torus $Z$.
Since $Z$ is a torus,
$\pi(\phi\inverse(X))$ is $Z$-cr (Lemma \ref{lem:linred=>gcr}).
Thus, by Lemma \ref{lem:bmr2.12}(i), $\phi\inverse(X)$ is $\tilde{G}$-cr if
and only if $X_i$ is $G_i$-cr for each $i$.
This proves the first assertion.
If $X$ is connected, then for each $i$,
$X_i$ is generated by commuting subgroups $X_i^0$
and $\pi_i({\rm ker}(\phi))$.
Since $\pi_i({\rm ker}(\phi))$ is linearly reductive,
$X_i$ is $G_i$-cr if and only if $X_i^0$ is $G_i$-cr.
For, if $X_i$ is $G$-cr, then $X_i^0$ is $G$-cr by Theorem \ref{thm:BMR3.10}.
Conversely, if $X_i^0$ is $G$-cr, then $\pi_i({\rm ker}(\phi))$,
being linearly reductive, is $C_G(X_i^0)$-cr by Lemma \ref{lem:linred=>gcr},
so $X_i$ is $G$-cr by Proposition \ref{prop:AB<=>CGA}.
The second assertion of the lemma now follows from the first.
\end{proof}

\begin{lem}
\label{lem:isogeny}
  Let $f\colon G_1\ra G_2$ be an isogeny of connected reductive  
groups.  Then Theorem \ref{thm:commuting} holds for $G_1$ if and only  
if it holds for $G_2$.
\end{lem}

\begin{proof}
  Suppose Theorem \ref{thm:commuting} holds for $G_2$.  Let $A_1$,  
$B_1$ be connected commuting $G_1$-cr subgroups of $G_1$.  Then  
$A_2:=f(A_1)$ and $B_2:=f(B_1)$ are connected commuting $G_2$-cr  
subgroups of $G_2$ by Lemma \ref{lem:bmr2.12}(ii), so our hypothesis  
on $G_2$ implies that $A_2B_2$ is $G_2$-cr.  Since $f(A_1B_1)=A_2B_2$,  
Lemma \ref{lem:bmr2.12}(ii) implies that $A_1B_1$ is $G_1$-cr.

  Conversely, suppose Theorem \ref{thm:commuting} holds for $G_1$.   
Let $A_2$, $B_2$ be connected commuting $G_2$-cr subgroups of $G_2$.   
Then $A_1:=f^{-1}(A_2)^0$ and $B_1:=f^{-1}(B_2)^0$ are connected  
$G_1$-cr subgroups of $G_1$, by Lemma \ref{lem:bmr2.12}(ii).  Since  
$\ker f\subseteq Z(G_1)$ is linearly reductive, $A_1$ and $B_1$  
satisfy the hypotheses of Lemma \ref{lem:linredcomm}.  Hence $A_1B_1$  
is $G_1$-cr, by Lemma \ref{lem:linredcomm} and our hypothesis on  
$G_1$.  Lemma \ref{lem:bmr2.12}(ii) now implies that $A_2B_2=  
f(A_1B_1)$ is $G_2$-cr.
\end{proof}

\begin{lem}
\label{lem:simple}
  If Theorem \ref{thm:commuting} holds for each
simple factor $G_i$ of $G$, then it holds for $G$.
\end{lem}

\begin{proof}
Let $G$, $A$ and $B$ be as in the statement of Theorem \ref{thm:commuting}.
Multiplication gives an isogeny from $G_1\times \cdots \times  
G_r\times Z$ onto $G$,
where $Z=Z(G)^0$.
By Lemma \ref{lem:isogeny}, we may assume that
$G = G_1 \times \cdots \times G_r \times Z$.
Let $\pi\colon G\ra Z$ and $\pi_i\colon G\ra G_i$ be the projection maps.
Now $A$ (resp.\ $B$, $AB$) is $G$-cr if and only
if $\pi_i(A)$ (resp.\ $\pi_i(B)$, $\pi_i(AB)$)
is $G_i$-cr for each $i$, by Lemma \ref{lem:bmr2.12}(i)
(note that $\pi(A)$, $\pi(B)$ and $\pi(AB)$ are
automatically $Z$-cr as $Z$ is a torus).
But $\pi_i(AB) = \pi_i(A)\pi_i(B)$ and
$[\pi_i(A),\pi_i(B)] = \{1\}$, so Theorem \ref{thm:commuting} holds for $G$
if it holds for each $G_i$.
\end{proof}

We do not have a uniform proof of Theorem \ref{thm:commuting}; in fact
we proceed by a series of case-by-case checks.  First we consider
the classical groups.

\subsection{Classical Groups}
If $G$ is classical, then we obtain slightly stronger results
(see Remark \ref{rem:nonconn}).
First we consider the case $G=\GL(V)$.
We believe that the following result is a standard fact in
representation theory,
but we have not been able to find a proof in the literature.
The special case of Lemma \ref{lem:linearthm} when
$A$ and $B$ are connected reductive subgroups of
$\GL(V)$ is proved in \cite[Lem.\ 41]{mcninch1} using facts from the
representation theory of reductive groups.
We are grateful to R.\ Tange for providing the argument given below.

\begin{lem}
\label{lem:linearthm}
Let $V$ be a finite dimensional vector space over $k$.
Suppose $A$ and $B$ are commuting subgroups of $\GL(V)$.
Then $V$ is semisimple for the product $AB$ if and only if $V$ is semisimple
for $A$ and $B$.
\end{lem}

\begin{proof}
Suppose that $V$ is semisimple for $A$ and $B$.
Let $C$, $D$, and $E$ be the $k$-subalgebras of ${\rm End}(V)$
spanned by $A$, $B$, and $AB$, respectively.  Since $C$ and $D$
act faithfully and semisimply on $V$, $C$ and $D$ are
semisimple $k$-algebras (cf. \cite[Ch.\ XVII, Prop.\ 4.7]{Lang}).
Since $k$ is algebraically closed, and thus perfect,
$C\otimes_k D$ is semisimple, by
\cite[\S 7.6 Cor.\ 4]{bourbaki} (or \cite[Ch.\ XVII, Thm.\ 6.4]{Lang}).
We have an epimorphism from $C\otimes_k D$ to $E$ given by
$c\otimes d\mapsto cd$.  It follows that $E$ is also semisimple
(see \cite[Ch.\ XVII, Prop.\ 2.2 and \S 4]{Lang}), so $E$ acts
semisimply on $V$ (\cite[Ch.\ XVII, Prop.\ 4.7]{Lang}).
Thus $V$ is semisimple for $AB$, as required.
The other implication follows from Clifford's Theorem.
\end{proof}

\begin{thm}
\label{thm:commutingclassical}
Theorem \ref{thm:commuting} holds for $G$ a simple group of classical type.
\end{thm}

\begin{proof}
Let $A$, $B$ be commuting subgroups of $G$.  Let $G_1$ be the classical
group with the same Dynkin type as $G$: so $G_1$ is either $\SL(V)$,
$\SP(V)$ or $\SO(V)$.  Let $\widetilde{G}$ be the simply connected
cover of $G$.  We have canonical projections $\widetilde{G}\ra G$ and
$\widetilde{G}\ra G_1$.  By Lemma \ref{lem:isogeny}, we  
can assume that $G=G_1$.

If $G = \SL(V)$, then $G$ is normal in $\GL(V)$, so the result follows
from Proposition \ref{prop:10} and Lemma \ref{lem:linearthm}.
The other two possibilities $G=\SP(V)$ and $G=\SO(V)$
arise as the connected centralizer of an involution
acting on $\GL(V)$.
By hypothesis, $p \neq 2$, so the group of automorphisms generated by  
this involution
is linearly reductive
and Proposition \ref{prop:linred} applies; cf.\ \cite[Ex.\ 3.23]{BMR}.
Now $A$ and $B$ are $G$-cr if and only if $A$ and $B$ are
$\GL(V)$-cr, which happens if and only if $V$ is semisimple as an
$A$- and a $B$-module.
By Lemma \ref{lem:linearthm}, this happens if and only if $V$
is semisimple for $AB$,
which occurs if and only if $AB$ is $\GL(V)$-cr;
by Proposition \ref{prop:linred} again, this happens
if and only if $AB$ is $G$-cr.
\end{proof}

\begin{rem}
\label{rem:nonconn}
We did not need to assume that $A$ and $B$ were connected in Lemma  
\ref{lem:linearthm}.
The proofs of Lemmas \ref{lem:isogeny} and \ref{lem:simple} and  
Theorem \ref{thm:commutingclassical} also go through for non-connected  
$A$ and $B$, so Theorem \ref{thm:commuting}
holds for non-connected $A$ and $B$ as well when all of the simple  
factors of $G$ are classical.
\end{rem}

\subsection{Exceptional Groups}

Recall Theorem \ref{thm:ls1}.
Examination of the possible values for $N(X,G)$ in Table \ref{t:nxg}
shows that $N(X,G) \leq 7$ always.
Thus, if $p>7$ (recall our convention that $p=\infty$ in  
characteristic zero!),
Theorem \ref{thm:commuting} holds for simple exceptional groups
simply because $AB$ is connected reductive.
The remainder of this subsection is devoted to improving the bound on  
$p$; we show that $p>3$ will do.

We now prove Theorem \ref{thm:commuting} via a series of lemmas
which exhaust all further possibilities.  At various points we
use inductive arguments involving Levi subgroups of groups of
type $E_6$, $E_7$ and $E_8$; we are able to leave the $G_2$ and
$F_4$ cases until last (Lemma \ref{lem:FG}) because groups of
type $G_2$ and $F_4$ cannot arise as simple factors of these Levi
subgroups.

\begin{lem}
\label{lem:E6E7rank4}
  Suppose $p>3$ and
$X$ is a connected reductive group which has simple
factors of rank at most $4$.
Suppose further that $X$ has no simple factor of type $A_4$ or $C_4$.
Then, if $Y$ is a connected reductive subgroup of $X$,
$Y$ is $X$-completely reducible.
\end{lem}

\begin{proof}
Using Lemma \ref{lem:simplefactors}, we reduce to the case
when $X$ is simple.
Thus we can list the possible types for $X$:
$A_1$, $A_2$, $B_2$, $G_2$, $A_3$, $B_3$, $C_3$, $B_4$, $D_4$, $F_4$.

First suppose $X$ has one of the classical types in this list.
Examining \cite[Table 8.4]{liebeckseitz}, we see that all these
types arise as subgroups of $F_4$ when $p>3$.
By Table \ref{t:nxg},
$N(Y,F_4) \leq 3 < p$ for any connected reductive subgroup $Y$ of $F_4$.
Since $X$ has classical type, we can invoke Theorem \ref{thm:ls3.8}
to conclude that $Y$ is $X$-cr.

Finally, if $X$ has type $G_2$ or $F_4$, then $N(Y,X) \leq 3 < p$,
by Table \ref{t:nxg}, so $Y$ is $X$-cr, by
a direct application of Theorem \ref{thm:ls1}.
\end{proof}

Now we introduce some more notation
to make the exposition easier.
Given a $G$-completely reducible subgroup $A$
of the simple exceptional group $G$,
let $H_A := C_G(A)^0$.
If $A$ and $B$ are connected, commuting $G$-cr subgroups of $G$, then
to prove that $AB$ is $G$-cr, it suffices to show that $B$ is $H_A$-cr,
by Propositions \ref{prop:finiteindex} and \ref{prop:AB<=>CGA} (note  
that $C_G(A)$ need not be connected).
The next lemma allows us to proceed through each remaining
simple group in turn.

\begin{lem}
\label{lem:finite}
Suppose $G$ is simple and Theorem \ref{thm:commuting}
holds for all simple factors of all proper Levi subgroups of $G$.
Then we may assume when proving Theorem \ref{thm:commuting} that $A  
\cap H_A$ is finite.
\end{lem}

\begin{proof}
If $A\cap H_A$ is an infinite group, then $A\cap H_A$ contains a non-trivial
(hence non-central) torus $S$ of $G$.
In this case, $AB \subseteq AH_A \subseteq C_G(S)$,
and $L:=C_G(S)$ is a proper Levi
subgroup of $G$.
By Corollary \ref{cor:levi}, $A$ and $B$ are commuting $L$-cr  
subgroups of $L$;
thus, by Lemma \ref{lem:simple} and the hypothesis, $AB$ is $L$-cr,
and hence is $G$-cr by Corollary \ref{cor:levi}.
\end{proof}

\begin{lem}
\label{lem:E6}
Theorem \ref{thm:commuting} holds for $G$ simple of type $E_6$.
\end{lem}

\begin{proof}
Let $A$ and $B$ be commuting $G$-cr subgroups of $G$.
We show that $B$ is $H_A$-cr.
Since any simple factor of a proper Levi subgroup of $G$ is of classical
type, we may assume $A \cap H_A$ is finite,
by Lemma \ref{lem:finite} and Theorem \ref{thm:commutingclassical}.
In particular, $A$ is not a non-trivial torus and $\rk AH_A = \rk A +  
\rk H_A$.
If $\rk AH_A = \rk G$, then $AH_A$ is a regular reductive subgroup
of $G$, and hence $H_A$ is a connected regular reductive subgroup of $G$.
Since $p$ is good for $G$ and $B$ is $G$-cr, Proposition \ref{prop:regular}
shows that $B$ is $H_A$-cr, as required.
Further, if $\rk B = \rk H_A$, then $B$ is regular in
$H_A$ and hence $H_A$-cr, by
Lemma \ref{lem:reg=>gcr}.  Also, the result is trivial if $A$ or $B$  
is trivial,
so we assume this is not the case.
We are therefore left to consider the cases where
$1 \leq \rk B < \rk H_A < \rk AH_A< \rk G = 6$.
Thus $\rk H_A \leq 4$.

Now Lemma \ref{lem:E6E7rank4} covers all
these cases except $H_A = A_4$ or $C_4$.
An examination of \cite[Table 8.3]{liebeckseitz} shows that if $H_A =  
A_4$, then
$H_A$ is regular, so that $B$ is $H_A$-cr by Proposition \ref{prop:regular}.
On the other hand, if $H_A = C_4$, then $A \subseteq C_G(H_A)^0 =  
\{1\}$, a contradiction,
which shows that this case cannot arise.
\end{proof}

\begin{lem}
\label{lem:E7}
Theorem \ref{thm:commuting} holds for $G$ simple of type $E_7$.
\end{lem}

\begin{proof}
Any simple factor of a proper Levi subgroup of $G$ is either of  
classical type or of type $E_6$.  By  Theorem  
\ref{thm:commutingclassical}
and Lemmas \ref{lem:finite} and \ref{lem:E6}, we may assume that
$A\cap H_A$ is finite.
Repeating the rank argument in the proof of Lemma \ref{lem:E6},
this means we only need to consider
cases where $1\leq \rk B < \rk H_A < \rk AH_A < \rk G = 7$ and $A$ is  
not a torus.

We now look at the
possibilities for $H_A$ with $\rk H_A \leq 5$.  For each one, we prove  
that $AB$ is $G$-cr, either directly or by proving that $B$ is  
$H_A$-cr.  Again, Lemma \ref{lem:E6E7rank4} covers most of the cases;
we are left to consider the possibility that
$$
H_A = A_4, C_4, A_4T_1, C_4T_1 , A_4A_1, C_4A_1, A_5, B_5, C_5, D_5,
$$
where $T_1$ denotes a $1$-dimensional torus.
We deal with these cases by examining \cite[Table 8.2]{liebeckseitz}.
\begin{itemize}
\item If $H_A=A_4, A_5$ or $D_5$, then $H_A$ is regular, so $B$ is  
$H_A$-cr by Proposition \ref{prop:regular}.
\item If $H_A = C_4$, then $A \subseteq C_G(H_A) = T_1$
is a torus.  But this is impossible, because we assume $A$ is not a torus.

This also shows that the case $H_A = C_4A_1$ cannot arise,
as there is no $A_1$ subgroup centralizing a $C_4$.
\item If $H_A = A_4T_1$ or $C_4T_1$, then $AB \subseteq C_G(T_1)$,
which is a proper Levi subgroup of $G$,
so $AB$ is $G$-cr by Corollary \ref{cor:levi}, Lemma \ref{lem:simple},  
Theorem \ref{thm:commutingclassical} and Lemma \ref{lem:E6}.
\item If $H_A = A_4A_1$, then since $\rk A + \rk H_A = \rk AH_A < \rk G=7$,
we have $\rk A = 1$. Then $A = A_1$ and, since $A\cap H_A$ is finite,
there is a subgroup of type $A_1A_1$
in $C_G(A_4)^0$.  But $C_G(A_4)^0 = A_2T_1$, which does not contain an  
$A_1A_1$ subgroup.
Thus this case cannot arise.
\item If $H_A = B_5$, then $H_A \subset D_6$ (cf.\ Remark  
\ref{rem:useoftables}), which is regular in $G$.
Moreover, the $B_5$ is the connected centralizer of the graph automorphism
of order $2$ for $D_6$.
Now $B$ is $D_6$-cr, by Proposition \ref{prop:regular} (note that  
$p>2$), and hence
$H_A$-cr, by Proposition \ref{prop:linred}.
\item There is no subgroup of type $C_5$ in $E_7$.
\end{itemize}
This completes the proof of the lemma.
\end{proof}

Now we need to deal with the case when $G$ has 
type $E_8$.
This is more involved because we actually allow $p$ to be 5,
which is a bad prime for $E_8$.

\begin{lem}
\label{lem:E8reg}
Suppose $p>3$.  Let $G$ be simple of type $E_8$ and
let $H$ be a connected reductive regular subgroup of $G$
such that $H$ is not simple of rank 8.  Then for any subgroup $K$ of $H$,
$K$ is $G$-completely reducible if and only if $K$ is $H$-completely  
reducible.
\end{lem}

\begin{proof}
Let $H_1,\ldots,H_r$ be the simple factors of $H$ and let $K_i$ be the
projection of $K$ to each $H_i$.  By Lemma \ref{lem:simplefactors}, it  
is enough
to prove the result for each $i$ with $H$ replaced by $H_i$ and $K$
replaced by $K_i$.  Now $H_i$ has semisimple rank at most 7 by hypothesis,
so $H_i$ is a regular subgroup of a proper Levi subgroup $L$ of $G$.
Since $p>3$ and $L$ has no simple factors of type $E_8$, $p$ is good for $L$,
so the required result follows from Proposition \ref{prop:regular}
and Corollary \ref{cor:levi}.
\end{proof}

\begin{lem}
\label{lem:E8cr=>Xcr}
Let $G$ be simple of type $E_8$ and let $X$ be a simple subgroup of  
$G$ such that $\rk X\leq 6$,
$X$ is not of type $C_4$ and $X$ is not a non-regular subgroup
of type $A_4$.  Then for any connected reductive subgroup $Y$ of $X$,
if $Y$ is $G$-completely reducible, then $Y$ is $X$-completely reducible.
\end{lem}

\begin{proof}
By Lemma \ref{lem:E6E7rank4}, we need only consider the cases when $X$
is regular and of type $A_4$ or $X$ has rank either 5 or 6.  We deal  
with these
cases by examining \cite[Table 8.1]{liebeckseitz}.  Let $Y$ be a
$G$-cr subgroup of $X$.
  \begin{itemize}
   \item If $X = B_5$, then $X \subset D_6$ (cf.\ Remark  
\ref{rem:useoftables}), which is regular in $G$.
Moreover, the $B_5$ subgroup of $D_6$ arises as the connected centralizer of
an involution of $D_6$.
By Lemma \ref{lem:E8reg}, $Y$ is $D_6$-cr,
so by Proposition \ref{prop:linred}, $Y$ is $X$-cr (note that $p>2$).
The same argument works for
$B_6 \subset D_7 \subset G$, where $D_7$ in $G$ is regular.
   \item If $X = A_5, A_6, D_5, D_6$ or $E_6$, then $X$ is regular in $G$, so
we are done by Lemma \ref{lem:E8reg}.
If $X = A_4$, then $X$ is regular in $G$
by hypothesis, so the same argument holds.
   \item There is no subgroup of type $C_5$ or $C_6$ in $E_8$.
  \end{itemize}
This completes the proof of the lemma.
\end{proof}

\begin{lem}
\label{lem:E8}
Theorem \ref{thm:commuting} holds for $G$ simple of type $E_8$.
\end{lem}

\begin{proof}
  Any simple factor of a proper Levi subgroup of $G$ is either of  
classical type or of type $E_6$ or $E_7$.  By  Theorem  
\ref{thm:commutingclassical}
and Lemmas \ref{lem:finite}, \ref{lem:E6} and \ref{lem:E7}, we may assume that
$A\cap H_A$ is finite.  We can assume that $A$ and $H_A$ have rank at  
least 1.  Hence if $AH_A$ is regular
in $G$, then, since $B$ is $G$-cr, Lemma \ref{lem:E8reg} implies that  
$B$ is $AH_A$-cr, so $B$ is $H_A$-cr by Proposition \ref{prop:10}.   
Thus we may assume that $\rk AH_A < \rk G$.
Repeating the rank argument in the proof of Lemma \ref{lem:E6},
this means we only need to consider
cases where $1\leq \rk B < \rk H_A < \rk AH_A < \rk G = 8$ and $A$ is  
not a torus.

We now look at the
possibilities for $H_A$ with $\rk H_A \leq 6$.  For each one, we prove  
that $AB$ is $G$-cr, either directly or by proving that $B$ is  
$H_A$-cr.  We use \cite[Table 8.1]{liebeckseitz} to deal with the
various cases.

  \begin{itemize}
   \item Suppose $H_A$ is not semisimple: say $S$
is a non-trivial central torus in $H_A$.
Then $AB$ is contained in $C_G(S)$, a proper Levi subgroup of $G$, so  
$AB$ is $G$-cr by Corollary \ref{cor:levi}, Lemma \ref{lem:simple},  
Theorem \ref{thm:commutingclassical}, Lemma \ref{lem:E6} and Lemma  
\ref{lem:E7}.
   \item Suppose $H_A$ has a non-regular $A_4$-factor.
Then $C_G(H_A)\subseteq C_G(A_4)$ is trivial; but this is impossible,
because $A$ is non-trivial.
   \item Suppose $H_A=C_4$.  There are two cases.  First, suppose $H_A$
is contained in a regular $E_6$.  The connected centralizers of $H_A$
and the $E_6$ are the same: this centralizer is of type $A_2$.  Thus $AH_A$
is contained in $M:=E_6A_2$, which has rank 8 and hence is regular in $G$.
By Lemma \ref{lem:E8reg}, $A$ and $B$ are $M$-cr, so $AB$ is $M$-cr by  
Lemma \ref{lem:simple}, Theorem \ref{thm:commutingclassical}
and Lemma \ref{lem:E6}.  Hence $AB$ is $G$-cr by Lemma \ref{lem:E8reg}.

   Second, suppose $H_A$ is contained in a regular $A_7$.  Then $H_A$
is the connected centralizer of an involution of the $A_7$ (cf.\  
Remark \ref{rem:useoftables}), so
Lemma \ref{lem:E8reg} and Proposition \ref{prop:linred} imply
that $B$ is $H_A$-cr.
   \item Suppose $H_A=C_4H$, where $H$ is connected and reductive but  
not a torus.
Then $H$ is a subgroup of $C_G(C_4)^0= A_2$ and $H$
is not a torus, so $C_G(H_A)^0= C_{A_2}(H)^0$ is a torus.  But $A$ is a
subgroup of $C_G(H_A)^0$ of semisimple rank at least one, a contradiction.
   \item Suppose $H_A$ is simple and not of type $A_4$ or type $C_4$.   
Since $H_A$ has rank at most 6, Lemma \ref{lem:E8cr=>Xcr} implies that  
$B$ is $H_A$-cr.
   \item Suppose every simple factor of $H_A$ has rank at most 4 and  
$H_A$ has no simple factors of type $A_4$ or $C_4$.  Lemma  
\ref{lem:E6E7rank4} implies that $B$ is $H_A$-cr.
   \item Otherwise, write $H_A=X_1X_2$, where $X_2$ is semisimple and  
has every simple factor of rank at most 3, and $X_1$ is semisimple and  
has every simple
factor of rank at least 4.
Then $X_1$ is simple, and we may assume that $X_1\neq C_4$ and $X_1$  
is not a non-regular $A_4$.
Note that $X_1\neq C_5,C_6$ since $G$ does not contain a $C_5$ or a  
$C_6$.  Suppose first that $X_1$ is regular and of type $A_4$, or is  
of type $A_5$ or $D_5$.  Then $AB$ is a subset of $X_1C_G(X_1)^0$,  
which is of type $A_4A_4$, $A_5A_1A_2$ or $D_5A_3$ respectively.  Thus  
$X_1C_G(X_1)^0$ is regular in $G$.  It now follows from Lemma  
\ref{lem:simple}, Theorem \ref{thm:commutingclassical} and Lemma  
\ref{lem:E8reg} that $AB$ is $G$-cr (compare the first $H_A=C_4$ case  
above).

Now suppose that $X_1$ is of type $A_6$, $B_6$, $D_6$ or $E_6$.  Then  
$H_A=X_1$ since $H_A$ has rank at most 6, so $B$ is $H_A$-cr by Lemma  
\ref{lem:E8cr=>Xcr}.

Finally, suppose that $X_1$ is of type $B_5$.  If $H_A=B_5$ then $B$  
is $H_A$-cr by Lemma \ref{lem:E8cr=>Xcr}.  Otherwise we must have  
$H_A=B_5A_1$, since $H_A$ has rank at most 6.  Our hypotheses on the  
rank of $AH_A$ now imply that $A=A_1$, so the $B_5$-factor is  
centralized by an $A_1A_1$.  Now $B_5$ sits inside a subgroup of $G$  
of type $D_6$, and this $D_6$ has a group of type $A_1A_1$ as its  
connected centralizer.  Thus we have two groups of type $A_1A_1$  
inside $C_G(B_5)^0=B_2$.  As $p>2$, these two subgroups must be  
$B_2$-conjugate.  Hence $AB$ is contained in $D_6A_1A_1$, which is  
regular.  It now follows from Lemma \ref{lem:simple}, Theorem  
\ref{thm:commutingclassical} and Lemma \ref{lem:E8reg} that $AB$ is  
$G$-cr (compare the first $H_A=C_4$ case above).
  \end{itemize}
This exhausts all the possibilities, so the lemma is proved.
\end{proof}

\begin{lem}
\label{lem:FG}
Theorem \ref{thm:commuting} holds for $G$ simple of type $G_2$ or $F_4$.
\end{lem}

\begin{proof}
Suppose $G$ is simple of type $G_2$. Then $G$ arises
as the connected centralizer of $D_4$ under the
triality graph automorphism.
Since $p>3$, this automorphism generates a linearly reductive group.
Thus a subgroup of $G$ is $G$-cr if and only if it is $D_4$-cr,
by Proposition \ref{prop:linred}.
But Theorem \ref{thm:commuting} holds for $D_4$ by Theorem  
\ref{thm:commutingclassical},
so it holds for $G$.

Now suppose $G$ is simple of type $F_4$.
In this case $G$ arises as the connected centralizer of an involution  
of $E_6$.
Since $p>2$, this automorphism generates a linearly reductive group.   
Thus a subgroup of $G$ is $G$-cr if and only if it is $E_6$-cr, by
Proposition \ref{prop:linred}.
But Theorem \ref{thm:commuting} holds for $E_6$ by Lemma \ref{lem:E6},
so it holds for $G$.
\end{proof}

Theorem \ref{thm:commuting} now follows from Lemma \ref{lem:simple},
together with
Theorem \ref{thm:commutingclassical} and Lemmas \ref{lem:E6},  
\ref{lem:E7}, \ref{lem:E8} and \ref{lem:FG}.

Combining Theorems \ref{thm:BMR3.10} and \ref{thm:commuting} yields
the following.

\begin{cor}\label{cor:commutingpgood}
Suppose that $G$ is connected and that $p$
is good for $G$ or $p>3$.
If $A$ and $B$ are commuting connected reductive subgroups of $G$,
then $AB$ is $G$-completely reducible if and only if $A$ and $B$ are
$G$-completely reducible.
\end{cor}

We provide a reformulation of Theorem \ref{thm:commuting} which shows more
clearly its relation to Theorem \ref{thm:BMR3.10}.

\begin{cor}
Suppose that $G$ is connected and that $p$ is good for $G$ or $p>3$.
Suppose $N \subseteq H$ are connected subgroups of $G$
such that $N$ is normal in $H$.
If there exists a connected $G$-completely reducible subgroup $M$ of $C_G(N)$
such that $H = MN$, then $H$ is $G$-completely reducible if and only if
$N$ is $G$-completely reducible.
\end{cor}

\begin{proof}
If $H$ is $G$-cr, then $N$ is $G$-cr by Theorem \ref{thm:BMR3.10}.
On the other hand, if $N$ is $G$-cr, then Theorem \ref{thm:commuting} applies,
with $A = N$, $B = M$ and $AB = H$.
\end{proof}

\begin{rem}
\label{rem:counter}
One consequence of Theorem \ref{thm:commuting} is
that the converse of Corollary \ref{cor:CGAcr=>Gcr}
is true for connected groups if $p$ is good or $p>3$;
this follows from Proposition \ref{prop:AB<=>CGA}.
However, we give examples in Section \ref{sec:examples}
which show that the converse to
Corollary \ref{cor:CGAcr=>Gcr}, and hence Theorem \ref{thm:commuting},
is false without the restriction on the characteristic.
\end{rem}

\begin{rem}
  By Lemma \ref{lem:linredcomm}, Theorem \ref{thm:commuting} holds  
under the weaker hypotheses that $A$ and $B$ are $G$-cr, $A$  
normalizes $B$, $[A,B]$ is linearly reductive and $[A,B]$ centralizes  
$A$.
\end{rem}

\section{Counterexamples and Extensions}
\label{sec:examples}

As promised in Remark \ref{rem:counter},
in this section we provide
examples which show that Theorem \ref{thm:commuting}
fails in general without the hypotheses of connectedness
and good characteristic.
We also give an extension to Theorem \ref{thm:commuting} which shows that
one can remove the connectedness assumption at least in some cases.
Our first example shows that, even in good characteristic,
Theorem \ref{thm:commuting} fails for disconnected groups.

\begin{exmp}
\label{exmp:nonconn}
Suppose $p=2$ and $m\geq 4$ is even.
Define $\phi\in \Aut(\GL_{2m})$ by $\phi(g)= J(g^t)^{-1}J^{-1}$,
where $g^t$ denotes the matrix transpose of $g$, 
$J = \left(\begin{array}{cc}
0 & I_m \\ I_m & 0 \end{array}\right)$ and $I_m$ is the $m\times m$  
identity matrix.
Set $A = \langle \phi\rangle$ and $G = A\ltimes \GL_{2m}$.
Let $B = \SP_m$ embedded diagonally in the maximal rank subgroup
$\SP_m\times \SP_m$ of $\SP_{2m}$,
and consider the canonical embedding of $\SP_{2m}$ in $G^0$.
We can identify $\SP_{2m}$ with $C_G(A)^0$.
By \cite[Ex.\ 3.45]{BMR}, $B$ is $G^0$-cr but not $\SP_{2m}$-cr.

Observe that $\SP_{2m}$ is $G^0$-ir.
Thus if $\lambda\in Y(G^0) = Y(G)$ with
$A\SP_{2m}\subseteq P_\lambda$, then
$\lambda$ belongs to $Y(Z(G^0))$ and $R_u(P_\lambda)=\{1\}$, so
$A\SP_{2m}\subseteq P_\lambda= L_\lambda$.
This shows that $A\SP_{2m}$ is $G$-cr.
Thus the normal subgroup $A$ of $A\SP_{2m}$ is $G$-cr (Theorem  
\ref{thm:BMR3.10}).
By Proposition \ref{prop:finiteindex}, $B$ is $G$-cr.
However, $B$ is not $\SP_{2m}$-cr, so Proposition  
\ref{prop:finiteindex} implies that $B$ is not $A\SP_{2m}$-cr.
We have $A\SP_{2m}=C_G(A)$, so $B$ is not $C_G(A)$-cr.

Thus we have commuting subgroups $A$ and $B$ which are $G$-cr, but such
that $B$ is not $C_G(A)$-cr.
Hence, by Proposition \ref{prop:AB<=>CGA}, $AB$ is not $G$-cr.
In particular, even though $p=2$ is good for $G^0 = \GL_{2m}$,
Theorem \ref{thm:commuting} and the converse to
Corollary \ref{cor:CGAcr=>Gcr} fail for these subgroups
of the disconnected group $G$.
\end{exmp}

\begin{rem}
We cannot have an example of this kind inside a connected group:
for a non-trivial unipotent subgroup $A$ can never be $G$-cr if $G$
is connected (see Remark \ref{rem:gcr=>red}).
\end{rem}

The following example, due to M. Liebeck,
shows that Theorem \ref{thm:commuting}, and hence the converse to
Corollary \ref{cor:CGAcr=>Gcr}, can also fail for connected groups
in bad characteristic.

\begin{exmp}
\label{exmp:liebeck}
Suppose $p=2$.
We show that there exist connected commuting
subgroups $A,B$ of $\SP_8$ such that $A$, $B$ are $\SP_8$-cr but $AB$ is not.

Let $A=B=\SL_2$ and let $V_A$, $V_B$
be the natural modules for $A$, $B$ respectively.
Choose symplectic forms $(\cdot, \cdot)_A$,
$(\cdot, \cdot)_B$ for $V_A$, $V_B$ respectively.
Then $\SL(V_A)=\SP(V_A)$ and $\SL(V_B)=\SP(V_B)$.
Set $W := V_A\otimes V_B$ with the symplectic form $(\cdot, \cdot)_W$ given by
$$
(u_1\otimes u_2,v_1\otimes v_2)_W := (u_1,v_1)_A(u_2,v_2)_B.
$$
Then $A$, $B$ and
$A\times B$ act on $W$ in the obvious way,
and these actions preserve $(\cdot, \cdot)_W$.
Below we shall be interested in $A$-stable subspaces of $W$.
For $v\in V_B$, set
$$
V_A\otimes v := \{u\otimes v \mid  u\in V_A\},
$$
a subspace of $W$.
Since $W$ is the $A$-module direct sum of the irreducible $A$-modules
$V_A\otimes v_1$ and $V_A\otimes v_2$, where $v_1$ and $v_2$
are any two linearly independent vectors in $V_B$, we see that $A$
acts completely reducibly on $W$ and the proper non-trivial $A$-stable
subspaces of $W$ are precisely the subspaces of the form $V_A\otimes v$
for some $v\in V_B$.
In particular, the $A$-stable subspaces of $W$ have dimension $0$, $2$ or $4$.

There exists an $A$-module isomorphism
$\phi_A: V_A \to V_A^*$ corresponding to the symplectic form
on $V_A$.
Define $\phi_B: V_B \to V_B^*$ analogously and identify
$W^*$ with $V_A^*\otimes V_B^*$
via the $(A \times B)$-module isomorphism $\psi:=\phi_A\otimes \phi_B$;
this is precisely the isomorphism corresponding to the symplectic form on $W$,
and gives rise to a symplectic form $(\cdot,\cdot)_{W^*}$ on $W^*$
given by
\[
(\psi(w_1),\psi(w_2))_{W^*} := (w_1,w_2)_W.
\]

Consider the $(A\times B)$-module $U := W\oplus W^*$
endowed with the direct sum symplectic form,
which we denote by $(\cdot,\cdot)_U$.
The $(A\times B)$-action preserves the symplectic structure, so we can  
regard $A$, $B$ and $A\times B$ as subgroups of $\SP(U)$.
If $M$ is an $A$-stable subspace of $U$ then we have a short exact  
sequence of $A$-modules
$$
\{0\} \ra M\cap W \ra M \ra \pi_2(M)\ra \{0\},
$$ where
$\pi_2:U\ra W^*$ is the canonical projection.
Since the second and fourth term have even dimension,
$M$ must have even dimension.
In particular, if $M$ is isotropic then $M$ has dimension $0$, $2$ or $4$.
We claim that if $M$ is any $4$-dimensional $A$-stable
isotropic subspace of $U$,
then there exists an $A$-stable isotropic complementary subspace.
To establish this, we observe that such a subspace $M$
must either be of the form
$V_A\otimes v \oplus V_A^*\otimes f$ for some $0\neq v\in V_A$
and $0 \neq f \in V_A^*$, or of the form
$M_\theta := \graph(\theta)$, where $\theta : W \ra W^*$
is an $A$-module isomorphism (note that $W$ and $W^*$ are not isotropic).
Any subspace of the first type is isotropic, so we can take
a complement to be $V_A\otimes v'\oplus V_A^*\otimes f'$,
where $v$ and $v'$ (resp.\ $f$ and $f'$) are linearly
independent in $V_A$ (resp.\ $V_A^*$).
Given $M$ of the form $M_\theta$, choose any $0\neq v \in V_B$.
The $A$-stable subspace $\theta(V_A\otimes v)$ of $W^*$ is of the form  
$V_A^*\otimes f$ for some $f\in V_B^*$.
Choose $g\in V_B^*$ such that $f$ and $g$ are linearly independent;
then $M_\theta$ and $V_A\otimes v\oplus V_A^*\otimes g$ intersect  
trivially, and the latter subspace is isotropic.
This proves the claim.

We now prove that $A$ is $\SP(U)$-cr;
the analogous result for $B$ follows by symmetry.
The parabolic subgroups of $\SP(U)$ are precisely the stabilizers of flags
$$
{\mathcal F} : \{0\} \subset M_1\subset M_2\subset \cdots \subset M_r\subset U
$$
of isotropic subspaces of $U$.
Moreover, a parabolic subgroup
$P$ is opposite to the stabilizer of $\mathcal{F}$
if and only if $P$ is the stabilizer of a flag
$\{0\} \subset N_1\subset N_2\subset \cdots \subset N_r\subset U$
of isotropic subspaces such that $U = M_i^\perp\oplus N_i$ for each  
$i$.  (Both of these facts can easily be established by considering  
the parabolic subgroups containing the standard maximal torus $S$  
described in \cite[Ch.\ V, \S 23.3]{borel}.)
Let ${\mathcal F}$ be a flag of $A$-stable isotropic subspaces of $U$.
There are only three possible types of flag to check.
If ${\mathcal F}$ has the form $\{0\} \subset M\subset U$ with $\dim M = 4$,
then we are done, by the previous claim (note that $M=M^\perp$).
Now suppose that ${\mathcal F}$ has the form $\{0\}\subset M_1\subset  
M_2\subset U$,
where $\dim M_1=2$ and $\dim M_2 = 4$.
By the previous claim, there exists a $4$-dimensional $A$-stable
isotropic subspace $N_2$ of $U$ such that $U = M_2^\perp\oplus N_2$.
Since $A$ acts completely reducibly on $U$,
there exists an $A$-stable complement $N_1$
to $M_1^\perp\cap N_2$ in $N_2$.
Then $N_1$ is an isotropic $A$-stable
complement to $M_1^\perp$ and $N_1\subset N_2$, as required.
Finally, suppose that
${\mathcal F}$ has the form $\{0\}\subset M\subset U$ with $\dim M = 2$.
It is easy to show, by listing the possibilities for
$M$ as in the previous claim,
that $M$ is contained in an $A$-stable $4$-dimensional
isotropic subspace of $U$, and we can use the argument
of the second case to prove that $M^\perp$ has an $A$-stable
isotropic complement in $U$.
This completes the proof that $A$ is $\SP(U)$-cr.

We finish by proving that the isotropic $(A\times B)$-stable subspace
$M_\psi=\{(w,\psi(w)) \mid w\in W\}$
does not admit an isotropic $(A\times B)$-stable complement,
which proves that $A\times B$ is not $\SP(U)$-cr; we repeat
the argument of \cite[Ex.\ 3.45]{BMR} for the convenience of the reader.
Suppose that $N$ is such a complement.
It follows from the discussion above on $A$-stable subspaces
of $U$ that $N$ must be of the form $M_\theta$
for some $(A\times B)$-module isomorphism $\theta:W\ra W^*$.
Since $A\times B$ acts irreducibly on $W$ and $W^*$,
we must have $\theta=a\psi$ for some $a\in k$, by Schur's Lemma.
For $u,v\in W$, we have
\begin{align*}
((u,\theta(u)), (v,\theta(v)))_U &= ((u,a\psi(u)), (v,a\psi(v)))_U \\
                                  &= (u,v)_W + a^2(\psi(u),\psi(v))_{W^*}\\
                                  &= (1+a^2)(u,v)_W.
\end{align*}
As $N$ is isotropic, this expression is identically $0$,
so we must have $a=1$.  But then $N=M$, a contradiction.
\end{exmp}

\begin{rem}
\label{rem:exceptionaleg}
We can endow each of $V_A$ and $V_B$
(and hence the spaces $V_A^*$, $V_B^*$, $W$, $W^*$ and $U$)
in Example \ref{exmp:liebeck}
with quadratic forms compatible with the given symplectic forms.
The actions of $A$, $B$ and $A\times B$ on the various spaces are
compatible with these quadratic forms,
so we can regard $A$, $B$ and $A\times B$ as subgroups of
$\SO(W)$, $\SO(U)$, etc.;
for example, the image of $A\times B$ in $\SP(W)$ is precisely $\SO(W)$.
Parabolic subgroups of $\SO(U)$ correspond to stabilizers of flags of
totally singular subspaces of $U$.
An argument similar to that of Example \ref{exmp:liebeck} shows that
$A$ and $B$ are $\SO(U)$-cr but $A\times B$ is not.

We can use this result to provide counterexamples for
exceptional groups as well.
For example, the group $\SO_8(k)$ has type $D_4$ and the exceptional group
of type $E_6$ has a Levi subgroup of type $D_4$.
Since $G$-complete reducibility is invariant under taking isogenies, we
can view $A$, $B$ and $A\times B$ as subgroups of $E_6$ in this way.
The subgroups $A$ and $B$ are $D_4$-cr, hence they are $E_6$-cr, by 
Proposition \ref{prop:10} and Corollary \ref{cor:levi};
however, the product $AB$ is not $D_4$-cr, so cannot be $E_6$-cr, again
by Proposition \ref{prop:10} and Corollary \ref{cor:levi}.
\end{rem}

\begin{exmp}
\label{exmp:liebecknonconn}
We can modify Example \ref{exmp:liebeck} to obtain a similar counterexample
involving finite subgroups rather than connected ones.
Take $k$ to be the algebraic closure of the field with two elements.
We replace $A$ (resp.\ $B$) with the finite subgroup $A(q)$ (resp.\ $B(q)$),
where $q$ is a power of $2$;
for $q$ sufficiently large, \cite[Lem.\ 2.10]{BMR} implies that
$A(q)$ and $B(q)$ are $\SP_8$-cr but $A(q)B(q)$ is not.
\end{exmp}

Examples \ref{exmp:liebeck} and \ref{exmp:liebecknonconn}
show that even if $A$, $B$ and $A^0B^0$ are $G$-cr,
$AB$ need not be $G$-cr:
passing to finite extensions does not preserve $G$-complete reducibility.
We can, however, identify one special case in which this works.

\begin{prop}
\label{prop:index}
Let $A$, $B$ be $G$-completely reducible subgroups of $G$
such that $A$ normalizes $B$, $A^0$ centralizes $B^0$ and
the index of $(AB)^0$ in $AB$ is coprime to $p$.
Then $AB$ is $G$-completely reducible if and only if $(AB)^0$
is $G$-completely reducible.
In particular, $AB$ is $G$-completely reducible if $p$ is good for $G^0$.
\end{prop}

\begin{proof}
Since $AB/(AB)^0$ is linearly reductive, the first assertion follows from
Corollary \ref{cor:linred}(ii).
Clearly, $(AB)^0= A^0B^0$.
Since $A^0$, $B^0$ are $G^0$-cr by Theorem \ref{thm:BMR3.10}
and Proposition \ref{prop:finiteindex},
the second assertion follows from Theorem \ref{thm:commuting} and  
Proposition \ref{prop:finiteindex}.
\end{proof}

Note that Proposition \ref{prop:index} is consistent with
Example \ref{exmp:nonconn},
for in that case $p=2$ and the index of $(AB)^0$ in $AB$ is also $2$,
whence $AB/(AB)^0$ is not linearly reductive.

In our final examples we return to the case that $N$ is a normal subgroup
of $H$ and $M$ is a subgroup of $H$ such that $H = MN$ (cf. Section  
\ref{sec:clifford}).
We show that even if $H$ is $G$-cr and
$M$ is a complement to $N$, $M$ need not be $G$-cr (see the discussion  
following
Question \ref{qn:MnormalisesN}).

\begin{exmp}
\label{exmp:S3}
Let $p = 2$, let $H$ be the symmetric group $S_3$ embedded
irreducibly in $G = \GL_2$ and let $N$ be the subgroup of $H$ of order $3$.
Any subgroup of $H$ of order $2$ is a complement to $N$, but such a subgroup,
being unipotent, cannot be $G$-completely reducible
(Remark \ref{rem:gcr=>red}), since $G$ is connected.
\end{exmp}

\begin{exmp}
\label{exmp:badconn}
  Suppose $p=2$ and $m\geq 4$ is even.  Let $H$ be the maximal rank  
subgroup $\SP_m\times \SP_m$ of $G:=\SP_{2m}$, let $N=\{1\}\times  
\SP_m$, $M'=\SP_m\times \{1\}$ and let $M$ be $\SP_m$ embedded  
diagonally in $H$.  Then $H$, $N$ and $M'$ are $G$-cr, $M$ and $M'$  
are both complements to $N$ in $H$, but $M$ is not $G$-cr (see  
\cite[Ex.\ 3.45]{BMR}).  It is easy to check that $H=N_G(N)^0$, and it  
follows from Proposition \ref{prop:finiteindex} and Lemma  
\ref{lem:bmr2.12}(i) that $M$ is $N_G(N)$-cr.  This shows that  
Corollary \ref{cor:CGAcr=>Gcr} is false if we replace $C_G(N)$ with  
$N_G(N)$.
\end{exmp}

%%%%%%%%%%%%%%%%%%%%%%%%%%%%%%%%%%%%%%%%%%%%%%%%%%%%%%%%%%%%%%%%%%%%%%
%%%%%%%%%%%%% Acknowledgments
%%%%%%%%%%%%%%%%%%%%%%%%%%%%%%%%%%%%%%%%%%%%%%%%%%%%%%%%%%%%%%%%%%%%%%

\bigskip
{\bf Acknowledgements}:
We would like to thank M.W.\ Liebeck, C.W.\ Parker, G.R.\ Robinson,
and R.\ Tange for helpful discussions.
We are also grateful to the referees for their helpful comments.

The authors acknowledge the financial support of EPSRC Grant EP/C542150/1
and Marsden Grant UOC0501.
Part of the research for this paper was carried out while the
authors were staying at the Mathematical Research Institute
Oberwolfach supported by the ``Research in Pairs'' programme.
Also part of this research was carried out during a visit by the second author
to the University of Southampton: he is grateful to the
members of the School of Mathematics for their hospitality.

\bigskip

%%%%%%%%%%%%%%%%%%%%%%%%%%%%%%%%%%%%%%%%%%%%%%%%%%%%%%%%%%%%%%%%%%%%%%
%%%%%%%%%%%%% bibliography
%%%%%%%%%%%%%%%%%%%%%%%%%%%%%%%%%%%%%%%%%%%%%%%%%%%%%%%%%%%%%%%%%%%%%%

\end{document}